\documentclass[12pt]{amsart}

\usepackage{amssymb}

\newcommand{\be}{\begin{equation}} 
\newcommand{\ee}{\end{equation}}
\newcommand{\beq}{\begin{eqnarray}}
\newcommand{\eeq}{\end{eqnarray}}
\newcommand{\bt}{\beta}
\newcommand{\bl}{\begin{lemma}}
\newcommand{\el}{\end{lemma}}
\newcommand{\bm}{\begin{pmatrix}}
\renewcommand{\em}{\end{pmatrix}}
\newcommand{\bml}{\begin{multline}}
\newcommand{\eml}{\end{multline}}
\newcommand{\ba}{\begin{array}}
\newcommand{\ea}{\end{array}}

\newcommand{\la}{\label}
\newcommand{\ci}{\cite}

\newcommand{\de}{\delta}

\newcommand{\al}{\alpha}

\newcommand{\Ga}{\Gamma}

\newcommand{\om}{\omega}

\newcommand{\lb}{\lambda}

\renewcommand{\th}{\theta}

\newcommand{\bi}{\bibitem}

\newfont{\msbm}{msbm10 scaled\magstep1}
\newfont{\msbms}{msbm7 scaled\magstep1} 

\newcommand{\bbr}{\mbox{$\mbox{\msbm R}$}}

\newcommand{\bbc}{\mbox{$\mbox{\msbm C}$}}

\newcommand{\bbz}{\mbox{$\mbox{\msbm Z}$}}

\newtheorem{theorem}{Theorem}
\newtheorem{lemma}{Lemma}

\theoremstyle{definition}

\newtheorem{example}[theorem]{Example}

\theoremstyle{remark}
\newtheorem{remark}[theorem]{Remark}

\begin{document}
\def\wt{\widetilde}
\title[Toeplitz determinants]{Toeplitz and Hankel determinants 
with singularities: announcement of results}
\author{P. Deift}
\address{Courant Institute of Mathematical Sciences, New York, NY 10003, USA}  
\author{A. Its}
\address{Department of Mathematical Sciences,
Indiana University -- Purdue University  Indianapolis,
Indianapolis, IN 46202-3216, USA}
%
\author{I. Krasovsky}
\address{Department of Mathematical Sciences,
Brunel University West London,
Uxbridge UB8 3PH, United Kingdom}
%


\begin{abstract}
We obtain asymptotics for Toeplitz, Hankel, and Toeplitz+Hankel  determinants
whose symbols possess Fisher-Hartwig singularities. 
Details of the proofs will be presented in another publication.
\end{abstract}

\maketitle

Let $f(z)$ be a complex-valued function integrable over the unit circle with Fourier coefficients
\[
f_j={1\over 2\pi}\int_0^{2\pi}f(e^{i\theta})e^{-i j\theta}d\theta,\qquad j=0,\pm1,\pm2,\dots
\]
We are interested in the $n$-dimensional Toeplitz determinant with symbol $f(z)$,
\be\la{TD}
D_n(f(z))=\det(f_{j-k})_{j,k=0}^{n-1}.
\ee

In this paper we present the asymptotics of $D_n(f(z))$ as $n\to\infty$ in the case when the symbol $f(e^{i\theta})$ has a fixed number of 
Fisher-Hartwig singularities \ci{FH,L}, i.e., 
when it has the following form on the unit circle:
\be\la{fFH}
f(z)=e^{V(z)} z^{\sum_{j=0}^m \bt_j} 
\prod_{j=0}^m  |z-z_j|^{2\al_j}g_{z_j,\bt_j}(z)z_j^{-\bt_j},\qquad z=e^{i\th},\qquad
\theta\in[0,2\pi),
\ee
for some $m=0,1,\dots$,
where
\begin{eqnarray}
&z_j=e^{i\th_j},\quad j=0,\dots,m,\qquad
0=\th_0<\th_1<\cdots<\th_m<2\pi;&\la{z}\\
&g_{z_j,\bt_j}(z)\equiv g_{\bt_j}(z)=
\begin{cases}
e^{i\pi\bt_j}& 0\le\arg z<\th_j\cr
e^{-i\pi\bt_j}& \th_j\le\arg z<2\pi
\end{cases},&\la{g}\\
&\Re\al_j>-1/2,\quad \bt_j\in\bbc,\quad j=0,\dots,m,&
\end{eqnarray}
and $V(e^{i\theta})$ is a sufficiently smooth function on the unit circle. 
Here the condition on $\al$ insures integrability.
Note that a single Fisher-Hartwig singularity at $z_j$ consists of 
a root-type singularity
\be\la{za}
|z-z_j|^{2\al_j}=\left|2\sin\frac{\th-\th_j}{2}\right|^{2\al_j}
\ee
and a jump $g_{\bt_j}(z)$. 
A point $z_j$, $j=1,\dots,m$ is included in (\ref{z})
if and only if either $\al_j\neq 0$ or $\bt_j\neq 0$ (or both); 
in contrast, the point $z_0=1$ is always included (note that 
$g_{\bt_0}(z)=e^{-i\pi\bt_0})$.
Observe that for each $j$, $z^{\beta_j} g_{\beta_j}(z)$ 
is continuous at $z=1$, and so for each $j$ 
each ``beta'' singularity produces a jump only at the point $z_j$.
The factors $z_j^{-\bt_j}$ are singled out to simplify comparisons 
with existing literature. Indeed, (\ref{fFH}) with the notation 
$b(\th)=e^{V(e^{i\theta})}$ is exactly the symbol considered 
in \ci{FH,B,B2,BE1,BEsym,BEnonsym,BE4,BS,BS2,BS3,ES,Ehr,W}. 
We write the symbol, however, in a form with
$z^{\sum_{j=0}^m \bt_j}$ factored out. The present way of writing 
is more natural for our analysis.

On the unit circle $V(z)$ is represented by its Fourier expansion:
\be\la{fourier}
V(z)=\sum_{k=-\infty}^\infty V_k z^k,\qquad 
V_k={1\over 2\pi}\int_0^{2\pi}V(e^{i\th})e^{-ki\th}d\th.
\ee
The canonical Wiener-Hopf factorization of $e^{V(z)}$ is
\be\la{WienH}
e^{V(z)}=b_+(z) e^{V_0} b_-(z),\qquad b_+(z)=e^{\sum_{k=1}^\infty V_k z^k},
\qquad b_-(z)=e^{\sum_{k=-\infty}^{-1} V_k z^k}.
\ee

First, we address the (essentially known) case when all $\Re\bt_j$ lie in a single 
half-closed interval of length 1,
namely $\Re\bt_j\in(q-1/2,q+1/2]$, $q\in\bbr$,
reproving results of Szeg\H o for $\al_j=\bt_j=0$, 
Widom \ci{W} for $\bt_j=0$, Basor \ci{B} for $\Re\bt_j=0$,
B\"ottcher and Silbermann \ci{BS} for $|\Re\al_j|<1/2$, $|\Re\bt_j|<1/2$,
Ehrhardt \ci{Ehr} for $|\Re\bt_j-\Re\bt_k|<1$,
and other results of these authors (see \ci{Ehr} for a review). 
Note that we write the asymptotics in a form 
that makes it clear which branch of the roots is to be used.
\begin{theorem}\la{asTop}
Let $f(e^{i\theta})$ be defined in (\ref{fFH}) and $\al_j\pm\bt_j\neq -1,-2,\dots$
for $j=0,1,\dots$. Then as $n\to\infty$,
\begin{multline}\la{asD}
D_n(f)=\exp\left[nV_0+\sum_{k=1}^\infty k V_k V_{-k}\right]
\prod_{j=0}^m b_+(z_j)^{-\al_j+\bt_j}b_-(z_j)^{-\al_j-\bt_j}\\
\times
n^{\sum_{j=0}^m(\al_j^2-\bt_j^2)}\prod_{0\le j<k\le m}
|z_j-z_k|^{2(\bt_j\bt_k-\al_j\al_k)}\left({z_k\over z_j e^{i\pi}}
\right)^{\al_j\bt_k-\al_k\bt_j}\\
\times
\prod_{j=0}^m\frac{G(1+\al_j+\bt_j) G(1+\al_j-\bt_j)}{G(1+2\al_j)}
\left(1+o(1)\right),\\
\mbox{if}\qquad
\Re\al_j>-{1\over 2},\qquad |\Re\bt_j-\Re\bt_k|<1, \qquad j,k=0,1,\dots,m,
\end{multline}
where
$G(x)$ is Barnes' $G$-function. 
\end{theorem}

\begin{remark}
In the case of  a single singularity,
i.e., when $m=0$ or $m=1$, $\al_0=\bt_0=0$,
the theorem implies that the asymptotics (\ref{asD}) hold for
\be
\Re\al_m>-{1\over 2},\qquad \bt_m\in\bbc,\qquad \al_m\pm\bt_m\neq -1,-2,\dots
\ee
In fact, if there is only one singularity and $V\equiv 0$, an explicit formula is known \ci{BS}
for $D_n(f)$ in terms of the G-functions.
\end{remark}

\begin{remark}\la{error}
Assume that the function $V(z)$ is analytic. Then the following
can be said about the remainder term. If all $\bt_j=0$, the error term 
$o(1)=O(n^{-1}\ln n)$. If there is only one singularity 
the error term is also $O(n^{-1}\ln n)$. In the general case,
the error term depends on the differences $\bt_j-\bt_k$.
For analytic $V(z)$, our methods would allow us to calculate several asymptotic 
terms rather than just the main one presented in 
(\ref{asD}) (and also in (\ref{asDgen}) below). 
\end{remark}

\begin{remark}
If all $\Re\bt_j\in(-1/2,1/2]$ or all $\Re\bt_j\in[-1/2,1/2)$,
the conditions $\al_j\pm\bt_j\neq -1,-2,\dots$ are satisfied automatically
as $\Re\al_j>-1/2$.
\end{remark}

\begin{remark}\la{degen1}
Since $G(-k)=0$, $k=0,1,\dots$, the formula (\ref{asD}) no longer 
represents the leading asymptotics
if $\al_j+\bt_j$ or $\al_j-\bt_j$ is a negative integer for some $j$.
A similar situation arises in Theorem \ref{BT} below if 
some representations in $\mathcal{M}$ are degenerate.
We do not address this case in the paper.
\end{remark}

As mentioned above, Theorem \ref{asTop} was proved by Ehrhardt. We give 
an independent proof of this result using a connection of $D_n(f)$ with 
the system of 
polynomials orthogonal with weight $f(z)$ (\ref{fFH}) on the unit circle.
First, we can show that all $D_k(f)\neq 0$, $k=k_0,k_0+1\dots$,
for some sufficiently large $k_0$.
Then the polynomials 
$\phi_k(z)=\chi_k z^k+\cdots$, $\widehat\phi_k(z)=\chi_k z^{k}+\cdots$ 
of degree $k$, $k=k_0,k_0+1,\dots$, satisfying 
\begin{multline}\la{or0}
{1\over 2\pi}\int_0^{2\pi}\phi_k(z)z^{-j}f(z)d\theta=\chi_k^{-1}\de_{jk},\qquad
{1\over 2\pi}\int_0^{2\pi}\widehat\phi_k(z^{-1})z^j f(z)d\theta=
\chi_k^{-1}\de_{jk},\\
z=e^{i\theta},\qquad j=0,1,\dots,k,
\end{multline}
exist and are given by the following expressions:
\be\la{ef1}
\phi_k(z)={1\over\sqrt{D_k D_{k+1}}}
\left| 
\begin{matrix}
f_{00}& f_{01}& \cdots & f_{0k}\cr
f_{10}& f_{11}& \cdots & f_{1k}\cr
\vdots & \vdots &  & \vdots \cr
f_{k-1\,0} & f_{k-1\,1} & \cdots & f_{k-1\,k} \cr
1& z& \cdots & z^k
\end{matrix}
\right|,
\ee
\be\la{ef2}
\widehat\phi_k(z^{-1})={1\over\sqrt{D_k D_{k+1}}}
\left| 
\begin{matrix}
f_{00}& f_{01}& \cdots & f_{0\,k-1}& 1\cr
f_{10}& f_{11}& \cdots & f_{1\,k-1}& z^{-1}\cr
\vdots & \vdots &  & \vdots & \vdots\cr
f_{k0} & f_{k1} & \cdots & f_{k\,k-1}& z^{-k}
\end{matrix}
\right|,
\ee
where
\[ 
f_{st}={1\over 2\pi}\int_0^{2\pi}f(z)z^{-(s-t)}d\theta,\quad 
s,t=0,1,\dots,k.
\]
We obviously have
\be\la{chiD}
\chi_k=\sqrt{D_k\over D_{k+1}}.
\ee
These polynomials 
satisfy a Riemann-Hilbert problem. We solve the problem asymptotically for large $n$ in case of the weight given by (\ref{fFH}) with analytic $V(z)$, thus obtaining the large $n$ asymptotics of the orthogonal polynomials. The main new feature of the
solution is a construction of the local parametrix at the points $z_j$ of Fisher-Hartwig singularities. This parametrix is given in terms of the confluent hypergeometric function.
A study of the asymptotic behaviour of the polynomials 
orthogonal on the unit circle was initiated by Szeg\H o \ci{Szego}. 
Riemann-Hilbert methods developed within the last 20 years allow us to
find asymptotics of orthogonal polynomials in all regions 
of the complex plane (see \ci{Dstrong} and many subsequent works 
by many authors). 
Such an analysis of the polynomials with an analytic weight on the unit circle
was carried out in \ci{MMS1}, and for the case of a weight with 
$\al_j$-singularities but without jumps, in \ci{MMS2}. 
We provide, therefore, a generalization of these results.
Here we present only the following statement we will need below
for the analysis of determinants.

\begin{theorem}\la{poly}
Let $f(e^{i\theta})$ be defined in (\ref{fFH}), 
$V(z)$ be analytic in a neighborhood of the unit circle,
and $\phi_k(z)=\chi_k z^k+\cdots$, $\widehat\phi_k(z)=\chi_k z^{k}+\cdots$ be the corresponding polynomials satisfying (\ref{or0}).
Assume that 
$|\Re\bt_j-\Re\bt_k|<1$, $\al_j\pm\bt_j\neq -1,-2,\dots$, $j,k=0,1,\dots, m$.
Let
\be\la{de}
\de=\max_{j,k} n^{2\Re(\bt_j-\bt_k-1)}.
\ee
Then 
as $n\to\infty$,
\begin{multline}\la{aschi}
\chi_{n-1}^2=\exp\left[ -\int_0^{2\pi}V(e^{i\th}){d\th\over2\pi}\right]\left(
1-{1\over n}\sum_{k=0}^m (\al_k^2-\bt_k^2)\right.\\
+\left.
\sum_{j=0}^m\sum_{k\neq j}{z_k\over z_j-z_k}
\left({z_j\over z_k}\right)^n n^{2(\bt_k-\bt_j-1)}{\nu_j\over\nu_k}
{\Ga(1+\al_j+\bt_j)\Ga(1+\al_k-\bt_k)\over
\Ga(\al_j-\bt_j)\Ga(\al_k+\bt_k)}\frac{b_+(z_j)b_-(z_k)}{b_-(z_j)b_+(z_k)}
\right.\\
\left.
+ O(\de^2)+O(\de/n)\right),
\end{multline}
where 
\be
\nu_j=\exp\left\{-i\pi\left(\sum_{p=0}^{j-1}\al_p-
\sum_{p=j+1}^m\al_p\right)\right\}
\prod_{p\neq j}\left({z_j\over z_p}\right)^{\al_p}
|z_j-z_p|^{2\bt_p}.
\ee
Under the same conditions,
\be\la{asphi}
\phi_n(0)=\chi_n \left(
\sum_{j=0}^m n^{-2\bt_j-1} z_j^n \nu_j
{\Gamma(1+\al_j+\bt_j)\over \Gamma(\al_j-\bt_j)}
\frac{b_+(z_j)}{b_-(z_j)}+
O\left(\left[\de+{1\over n}\right]\max_k{n^{-2\Re\bt_k}\over n}
\right)\right),
\ee
\be\la{ashatphi}
\widehat\phi_n(0)=\chi_n
\left(
\sum_{j=0}^m n^{2\bt_j-1} z_j^{-n} \nu_j^{-1}
{\Gamma(1+\al_j-\bt_j)\over \Gamma(\al_j+\bt_j)}
\frac{b_-(z_j)}{b_+(z_j)}+
O\left(\left[\de+{1\over n}\right]\max_k{n^{2\Re\bt_k}\over n}
\right)\right).
\ee
\end{theorem}

\begin{remark}\la{17}
The error terms here are uniform and differentiable in all $\al_j$, $\bt_j$
for $\bt_j$ in compact subsets of the strip $|\Re\bt_j-\Re\bt_k|<1$, 
for $\al_j$ in compact
subsets of the half-plane $\Re\al_j>-1/2$, and outside a neighborhood
of the sets $\al_j\pm\bt_j=-1,-2,\dots$. If $\al_j+\bt_j=0$ or $\al_j-\bt_j=0$
for some $j$, the corresponding terms in the above formulas vanish.
\end{remark}
\begin{remark}
Note that the terms with $n^{2(\bt_k-\bt_j-1)}$ in (\ref{aschi}) become 
larger in absolute value that the $1/n$ term for $\max_{j,k}\Re(\bt_j-\bt_k)>1/2$.
\end{remark}

Our proof of Theorem \ref{asTop} uses Theorem \ref{poly}, similar 
results for the asymptotics of the orthogonal polynomials and their Cauchy 
transforms at the points $z_j$, and a set of differential identities 
for the logarithm of $D_n$, in the spirit of \cite{D,IK,Kduke}.

Our next task is to extend the result for arbitrary $\bt_j\in\bbc$,
i.e. for the case when not all $\Re\bt_j$'s lie in a single interval of 
length less than 1. We know from examples
(see, e.g. \ci{BS,BT,Ehr}) that in general, the formula (\ref{asTop}) breaks down.
Obviously, the general case can be reduced to $\Re\bt_j\in(q-1/2,q+1/2]$ by adding
integers to $\bt_j$. Then, apart from a constant factor, the only change in $f(z)$ is
multiplication with $z^\ell$, $\ell\in\bbz$. However, as can be shown, 
the determinants $D_n(f(z))$ and $D_n(z^\ell f(z))$ are simply related.
For example, for $\ell=1,2,\dots$,
\be\la{213}
D_n(z^\ell f(z))={(-1)^{\ell n} F_n\over \prod_{j=1}^{\ell-1}j!}
D_n(f(z)),
\ee
where
\[
F_n=\left|
\begin{matrix}
\Phi_n(0)& \Phi_{n+1}(0) & \cdots & \Phi_{n+\ell-1}(0)\cr
{d\over dz}\Phi_n(0) & {d\over dz}\Phi_{n+1}(0)& \cdots & 
{d\over dz}\Phi_{n+\ell-1}(0)\cr
\vdots & \vdots &  & \vdots\cr
 {d^{\ell-1}\over dz^{\ell-1}}\Phi_n(0) & {d^{\ell-1}\over dz^{\ell-1}}\Phi_{n+1}(0)& \cdots &
{d^{\ell-1}\over dz^{\ell-1}}\Phi_{n+\ell-1}(0)
\end{matrix}
\right|,
\]
and $\Phi_k(z)=\phi_k(z)/\chi_k$.
Since the $\chi_k$, $\phi_k(0)$, 
$\widehat\phi_k(0)$ for large $k$ are given by Theorem \ref{poly},
and expressions for the derivatives can be found similarly,
it is easy to obtain the general asymptotic formula for
$D_n$. However, this formula is implicit in the sense that one still
needs to separate the main asymptotic term from the others: e.g.,
if the dimension $\ell$ of $F_n$ is larger than the
number of leading-order terms in (\ref{asphi}), the obvious 
candidate for the leading order in $F_n$ vanishes (this
is not the case in the simplest situation given by Theorem
\ref{BT1}). We outline below how we resolve this problem.

Following \ci{BT,Ehr}, define a so-called representation 
of a symbol. Namely, for $f(z)$ given by (\ref{fFH}) replace 
$\bt_j$ by $\bt_j+n_j$, $n_j\in\bbz$
if $z_j$ is a singularity (i.e., if either $\bt_j\neq 0$ or 
$\al_j\neq 0$ or both: we set $n_0=0$ if $z_0=1$ is not a singularity).
The integers $n_j$ are arbitrary subject to the condition 
$\sum_{j=0}^m n_j=0$.
In a slightly different notation from \ci{BT,Ehr}, we call the 
resulting function $f(z;n_0,\dots,n_m)$ a representation of $f(z)$.
(The original $f(z)$ is also a representation corresponding to 
$n_0=\cdots=n_m=0$.)
Obviously, all representations of $f(z)$ differ only by 
multiplicative constants. 
We have
\be
f(z)=\prod_{j=0}^m z_j^{n_j}\times f(z;n_0,\dots,n_m).
\ee
We are interested in the representations (characterized by 
$(n_j)_{j=0}^m$) of $f$ such that
$\sum_{j=0}^m (\Re\bt_j+n_j)^2$ is minimal. There is a finite number 
of such representations and we provide an algorithm for finding 
them explicitly (see Remark \ref{Alg}). We call the set of 
them $\mathcal{M}$. Furthermore, we call a representation 
degenerate if $\al_j+(\bt_j+n_j)$ or $\al_j-(\bt_j+n_j)$ is a negative 
integer for some $j$. We call $\mathcal{M}$ non-degenerate if it contains 
no degenerate representations.
We prove 

\begin{theorem}\la{BT}
Let $f(z)$ be given in (\ref{fFH}), $\Re\al_j>-1/2$, $\bt_j\in\bbc$, $j=0,1,\dots,m$.
Let $\mathcal{M}$ be non-degenerate.
Then, as $n\to\infty$,
\be\la{asDgen}
D_n(f)=\sum\left(\prod_{j=0}^m z_j^{n_j}\right)^n 
\mathcal{R}(f(z;n_0,\dots,n_m))(1+o(1)),
\ee
where the sum is over all representations in $\mathcal{M}$.
Each $\mathcal{R}(f(z;n_0,\dots,n_m))$ stands for the right-hand side 
of the formula (\ref{asD}), without the error term,
corresponding to $f(z;n_0,\dots,n_m)$.
\end{theorem}

\begin{remark}
This theorem was conjectured by Basor and Tracy \ci{BT}.
The case when the representation minimizing 
$\sum_{j=0}^m (\Re\bt_j+n_j)^2$ is unique, i.e. there is only 
one term in the sum (\ref{asDgen}), was proved by Ehrhardt \ci{Ehr}.
Note that this case happens if and only if there exist such $n_j$ that
$\Re\bt_j+n_j$
belong to a half-open interval of length 1 for all $j=0,\dots,m$:
see next Remark.
Thus, Theorem \ref{BT} in this case follows from Theorem 
\ref{asTop} applied to this representation.
\end{remark}

\begin{remark}\la{Alg}
The set $\mathcal{M}$ can be characterized as follows. Suppose 
that the seminorm $\|\bt\|\equiv\max_{j,k}|\Re\bt_j-\Re\bt_k|>1$. Then, writing
$\bt^{(1)}_s=\bt_s+1$,  $\bt^{(1)}_t=\bt_t-1$, and  
$\bt^{(1)}_j=\bt_j$ if $j\neq s,t$,
where $\bt_s$ is one of the beta-parameters with
$\Re\bt_s=\min_j\Re\bt_j$,  $\bt_t$  is one of the beta-parameters with
$\Re\bt_t=\max_j\Re\bt_j$, we see that
$\|\bt^{(1)}\|\le\|\bt\|$, and $f$ corresponding to $\bt^{(1)}$ is 
a representation. After a finite number, say $r$, of such transformations
we reduce an arbitrary set of $\bt_j$ to the situation
for which either $\|\bt^{(r)}\|<1$ or $\|\bt^{(r)}\|=1$.
Note that further transformations do not change the seminorm in the second 
case, while in the first case the seminorm oscillates periodically 
taking 2 values, $\|\bt^{(r)}\|$ and $2-\|\bt^{(r)}\|$.
Thus all the symbols of type (\ref{fFH}) belong to 2 distinct classes:
the first, for which $\|\bt^{(r)}\|<1$, and the second, for which
$\|\bt^{(r)}\|=1$.
For symbols of the first class, 
$\mathcal{M}$ has only one member with beta-parameters $\bt^{(r)}$. 
Indeed, writing $b_j=\Re\bt_j$,
if $-1/2<b_j^{(r)}-q\le 1/2$ for some
$q\in\bbr$ and all $j$, then for any $(k_j)_{j=0}^m$ such 
that $\sum_{j=0}^m k_j=0$ and not all $k_j$ are zero, we have
\be\la{ineq}
\sum_{j=0}^m (b_j^{(r)}+k_j)^2=\sum_{j=0}^m (b_j^{(r)})^2+
2\sum_{j=0}^m(b_j^{(r)}-q)k_j+
\sum_{j=0}^m k_j^2
>\sum_{j=0}^m (b_j^{(r)})^2+\sum_{j=0}^m k_j^2-|k_j|\ge
\sum_{j=0}^m (b_j^{(r)})^2,
\ee
where the first inequality is strict as at least one $k_j>0$.
For symbols of the second class, we can find $q\in\bbr$ such that
$-1/2\le b_j^{(r)}-q\le 1/2$ for all $j$. Equation (\ref{ineq}) in this
case holds with ``$>$'' sign replaced by ``$\ge$''. Clearly, there are several
representations in $\mathcal{M}$ in this case (they correspond to
the equalities in (\ref{ineq})) and 
adding $1$ to one of $\bt^{(r)}_s$ with $b^{(r)}_s=\min_j b^{(r)}_j=
q-1/2$
while subtracting $1$ from one of $\bt^{(r)}_t$ with
$b^{(r)}_t=\max_j b^{(r)}_j=q+1/2$
provides the way to find all of them.

A simple explicit sufficient, but obviously not necessary, condition 
for $\mathcal{M}$ to have only one member is that all $\Re\bt_j\mod 1$
be different.
\end{remark}

\begin{remark} 
The situation when all $\al_j\pm\bt_j$ are nonnegative integers, which was 
considered by B\"ottcher and Silbermann in \ci{BS2}, is a particular case 
of the above theorem.
\end{remark}

\begin{remark}\la{degen2}
The case when {\it all} the representations of $f$ are degenerate
(not only those in $\mathcal{M}$) was considered by Ehrhardt \ci{Ehr} who found
that in this case $D_n(f)=O(e^{nV_0}n^r)$, where $r$ is any real number. 
We can reproduce this result by our methods but do not present it here.
\end{remark}

We prove Theorem \ref{BT} in the following way. Consider the set 
$\bt^{(r)}_j$ constructed in Remark \ref{Alg}. We have to consider only 
the second class, i.e. $\|\bt^{(r)}\|=1$. We then have, relabelling
$\bt^{(r)}_j$ according to increasing real part,
\[
\Re\bt^{(r)}_1=\cdots=\Re\bt^{(r)}_p<\Re\bt^{(r)}_{p+1}\le\cdots
\le\Re\bt^{(r)}_{m'-\ell}<
\Re\bt^{(r)}_{m'-\ell+1}=\cdots=\Re\bt^{(r)}_{m'},
\]
for some $p,\ell>0$. Here $m'=m+1$ if $z=1$ is a singularity,
otherwise $m'=m$. Now consider the symbol (not a representation of $f$)
$\wt f$ of type (\ref{fFH})
with beta-parameters denoted by $\wt\bt$ and given by
$\wt\bt_j=\bt^{(r)}_j$ for $j=1,\dots, m'-\ell$, and
$\wt\bt_j=\bt^{(r)}_j-1$ for $j=m'-\ell+1,\dots, m'$.
It is easy to see that the original symbol $f$ has ${\ell\choose\ell+p}$
representations in $\mathcal{M}$ obtained by shifting any $\ell$
out of $\ell+p$ parameters
$\wt\bt_j$ with the smallest real part to the right by $1$.  
Thus, $f(z)=c z^{\ell}\wt f(z)$, where $c$ is a simple constant factor.
To find the asymptotics of $D_n(f)$, we now use (\ref{213}) with
$f$ replaced by $\wt f$. The l.h.s. is then $c^{-n}D_n(f)$.
In the r.h.s. we have a factor $D_n(\wt f)$ to which Theorem \ref{asTop}
is applicable since $\|\wt\bt\|<1$, and an $\ell\times\ell$ determinant $F_n$
involving the polynomials orthogonal w.r.t. $\wt f$ (for simplicity,
consider $V(z)$ analytic in a neighborhood of the unit circle).
It is a crucial fact
that the size $\ell$ of this determinant is less than the number
of terms, $\ell+p$, in the expansion of $\phi_n(0)/\chi_n$ of the same
largest order $O(n^{-2\Re\wt\bt_1-1})$ (see (\ref{asphi}) with $\bt_j$
replaced by $\wt\bt_j$). This fact enables us to extract the 
leading asymptotic contribution to $F_n$ (resolving the problem 
mentioned above). The asymptotics of $F_n$
and $D_n(\wt f)$ combine together and produce (\ref{asDgen}).

We will now discuss a simple particular case of Theorem \ref{BT}
and present a direct independent proof in this case.

\begin{theorem}[A particular case of Theorem \ref{BT}]\la{BT1}
 Let the symbol $f^{\pm}(z)$ be obtained from $f(z)$ (\ref{fFH})
by replacing one $\bt_{j_0}$ with $\bt_{j_0}\pm 1$ for some 
fixed $0\le j_0\le m$.
Let $\Re\al_j>-{1\over 2}$, $\Re\bt_j\in (-1/2,1/2]$, $j=0,1,\dots,m$.
Then
\be
D_n(f^{+}(z))=z_{j_0}^{-n} {\phi_n(0)\over\chi_n}D_n(f(z)),\qquad
D_n(f^{-}(z))=z_{j_0}^n {\widehat\phi_n(0)\over\chi_n}D_n(f(z)).
\ee
These formulas together with
(\ref{asphi},\ref{ashatphi},\ref{aschi},\ref{asD}) yield the following
asymptotic description of $D_n(f^{\pm})$.
Let there be more than one singular points $z_j$ and all $\al_j\pm\bt_j\neq 0$.
For $f^{+}(z)$, let $\bt_{j_p}$, $p=1,\dots,s$ 
 be such that they have the same real part which is strictly less
than the real parts of all the other $\bt_j$, i.e. 
$\Re\bt_{j_1}=\cdots=\Re\bt_{j_s}<\min_{j\neq j_1,\dots,j_s}\Re\bt_j$. 
For $f^{-}(z)$ let one $\bt_{j_p}$, $p=1,\dots,s$ be 
such that  $\Re\bt_{j_1}=\cdots=\Re\bt_{j_s}>\max_{j\neq j_1,\dots,j_s}
\Re\bt_j$. 
Then the asymptotics of $D_n(f^{\pm})$ are given by the following: 
\be\la{asf}
D_n(f^{+})=z_{j_0}^{-n}\sum_{p=1}^s z_{j_p}^n\mathcal{R}_{j_p,+}(1+o(1)),\qquad
D_n(f^{-})=z_{j_0}^{n}\sum_{p=1}^s z_{j_p}^{-n}\mathcal{R}_{j_p,-}(1+o(1)),
\ee
where $\mathcal{R}_{j,\pm}$ is the right-hand side of (\ref{asD}) 
(without the error term) 
in which $\bt_j$ is replaced by $\bt_j\pm 1$, respectively.
\end{theorem}

\begin{proof}
For simplicity, we present the proof only for $V(z)$ analytic in
a neighborhood of the unit circle. 
Consider the case of $f^{-}(z)$. It corresponds to one of the $\bt_j$ shifted 
inside the interval $(-3/2,-1/2]$. Since
\[
z^{\sum_{j=0}^m\bt_j-1}=z^{-1}z^{\sum_{j=0}^m\bt_j},\qquad g_{\bt_{j_0}-1}(z)=
-g_{\bt_{j_0}}(z),\qquad z_{j_0}^{-\bt_{j_0}+1}=z_{j_0}z_{j_0}^{-\bt_{j_0}},
\]
we see that
\[
f^{-}(z)=-z_{j_0}z^{-1}f(z).
\]
Therefore, using the identity (an analogue of (\ref{213}) for $\ell=-1$)
\be\la{DD-1}
D_n(z^{-1} f(z))=(-1)^n{\widehat\phi_n(0)\over \chi_n}D_n(f),
\ee
we obtain
\[
D_n(f^{-}(z))=(-z_{j_0})^n D_n(z^{-1}f(z))=
z_{j_0}^n {\widehat\phi_n(0)\over\chi_n}D_n(f(z)).
\]
If, for some $j_1$, $j_2$, $\dots$, $j_s$, we have that 
$\Re\bt_{j_1}=\cdots=\Re\bt_{j_s}>\max_{j\neq j_1,\dots,j_s}\Re\bt_j$, 
then we see from
(\ref{ashatphi}) that only the addends with $n^{2\bt_{j_1}-1}$, $\dots$,
$n^{2\bt_{j_s}-1}$ give  contributions
to the main asymptotic term of $D_n(f^{-}(z))$. Using the relation 
$G(1+x)=\Ga(x) G(x)$, we obtain the formula (\ref{asf}) for $D_n(f^{-}(z))$.
The case of $f^{+}(z)$ is similar.
\end{proof}

\begin{example}
In \ci{BT} Basor and Tracy noticed a simple example of a symbol of
type (\ref{fFH}) 
for which the asymptotics of the determinant can be computed directly,
but are very different from (\ref{asD}). Up to a constant, the symbol is
\be
\wt f(e^{i\th})=
\begin{cases}  
-i,& 0<\th<\pi\cr
i,& \pi<\th<2\pi
\end{cases}. 
\ee
We can represent $\wt f$ as a symbol with $\bt$-singularities 
$\bt_0=1/2$, $\bt_1=-1/2$ at the points $z_0=1$ and $z_1=-1$, respectively:
\be\la{tildef}
\wt f(z)=g_{1,1/2}(z)g_{-1,-1/2}(z)e^{i\pi/2}
\ee
We see that $\wt f(z)=f^{-}(z)$ and $j_0=1$.
Therefore by the first part of Theorem \ref{BT1}, we have
\[
D_n(\wt f(z))=(-1)^n{\widehat\phi_n(0)\over\chi_n}D_n(f(z)),
\]
where $\phi_n(z)$, $\chi_n$, $D_n(f(z))$ correspond to $f(z)$ given by
(\ref{fFH}) with $m=1$, $z_0=1$, $z_1=e^{i\pi}$, $\bt_0=\bt_1=1/2$, 
$\al_0=\al_1=0$. 

Substituting (\ref{ashatphi},\ref{asD})
into the above equation for $D_n(\wt f(z))$, we obtain
\be
D_n(\wt f(z))=\frac{1+(-1)^n}{2}\sqrt{2\over n} G(1/2)^2 G(3/2)^2(1+o(1)),
\ee
which is the answer found in \ci{BT}. 

Alternatively, noting that $s=2$, $j_1=j_0=1$ and $j_2=0$ and
using (\ref{asf}) we obtain
\[
D_n(\wt f(z))=(-1)^n((-1)^n\mathcal{R}_{1,-}+\mathcal{R}_{0,-}).
\]
Since $\mathcal{R}_{1,-}=\mathcal{R}_{0,-}=(2n)^{-1/2}G(1/2)^2 G(3/2)^2(1+o(1))$, we
obtain the same result.

As noted by Basor and Tracy, $\wt f(z)$
has a different representation of type (\ref{fFH}), namely,
with $\bt_0=-1/2$, $\bt_1=1/2$, and we can write
\be\la{f2}
\wt f(z)=-g_{1,-1/2}(z)g_{-1,1/2}(z)e^{-i\pi/2}.
\ee
This fact was the origin of their conjecture. In the notation of Theorem \ref{BT},
the symbol (\ref{tildef}) has the two representations minimizing 
$\sum_{j=0}^1 (\bt_j+n_j)^2$, one with $n_0=n_1=0$ and the other with $n_0=-1$, $n_1=1$.
\end{example}

Note that in the case $\sum_{j=0}^m\bt_j=0$ 
we can always assume that $\Re\bt_j\in[-1/2,1/2]$.
The beta-singularities then are just piece-wise constant (step-like) 
functions. 
This case is relevant for our next result, which is on Hankel determinants.

Let $w(x)$ be an integrable complex-valued function on the interval $[-1,1]$. 
Then the Hankel determinant with symbol $w(x)$ is
\be
D_n(w(x))=\det\left(\int_{-1}^1 x^{j+k}
  w(x)dx\right)_{j,k=0}^{n-1}.
\ee
Define $w(x)$ for a fixed $r=0,1,\dots$ as follows:
\begin{multline}\la{wFH}
w(x)=e^{U(x)}\prod_{j=0}^{r+1}|x-\lb_j|^{2\al_j}\om_j(x)\\
1=\lb_0>\lb_1>\cdots>\lb_{r+1}=-1,\qquad 
\om_j(x)=
\begin{cases}
e^{i\pi\bt_j}& \Re x<\lb_j\cr
e^{-i\pi\bt_j}& \Re x>\lb_j
\end{cases},\qquad
\Re\bt_j\in(-1/2,1/2],\\
\bt_0=\bt_{r+1}=0,\qquad \Re\al_j>-{1\over 2}, \qquad j=0,1,\dots,r+1.
\end{multline}
where $U(x)$ is a sufficiently smooth function on the interval $[-1,1]$.
Note that we set $\bt_0=\bt_{r+1}=0$ without loss of generality as the 
functions $\om_0$, $\om_{r+1}$ are just constants on $(-1,1)$. 

We prove
\begin{theorem}\la{asHankel}
Let $w(x)$ be defined in (\ref{wFH}). Then as $n\to\infty$,
\begin{multline}\la{asDH}
D_n(w)=
D_n(1)e^{\left[(n+\al_0+\al_{r+1})V_0-\al_0 V(1)-\al_{r+1}V(-1)+
{1\over 2}\sum_{k=1}^\infty k V_k^2\right]}\\
\times
\prod_{j=1}^r b_+(z_j)^{-\al_j-\bt_j}b_-(z_j)^{-\al_j+\bt_j}\times
e^{\left[2i(n+A)\sum_{j=1}^r\bt_j\arcsin\lb_j+
i\pi\sum_{0\le j<k\le r+1}(\al_j\bt_k-\al_k\bt_j)
\right]}\\
\times
4^{-\left(An+\al_0^2+\al_{r+1}^2+\sum_{0\le j<k\le r+1}\al_j\al_k+\sum_{j=1}^r\bt_j^2
\right)}
(2\pi)^{\al_0+\al_{r+1}}
n^{2(\al_0^2+\al_{r+1}^2)+\sum_{j=1}^r(\al_j^2-\bt_j^2)}\\
\times
\prod_{0\le j<k\le r+1}|\lb_j-\lb_k|^{-2(\al_j\al_k+\bt_j\bt_k)}
\left|\lb_j\lb_k-1+\sqrt{(1-\lb_j^2)(1-\lb_k^2)}\right|^{2\bt_j\bt_k}\\
\times
{1\over G(1+2\al_0)G(1+2\al_{r+1})}
\prod_{j=1}^r(1-\lb_j^2)^{-(\al_j^2+\bt_j^2)/2}
\frac{G(1+\al_j+\bt_j) G(1+\al_j-\bt_j)}{G(1+2\al_j)}
\left(1+o(1)\right),\\
A=\sum_{k=0}^{r+1}\al_k,\qquad
\Re\al_j>-{1\over 2},\qquad \Re\bt_j\in\left(-{1\over 2},{1\over2}\right),
\qquad j=0,1,\dots,r+1,\qquad \bt_0=\bt_{r+1}=0,
\end{multline}
where $V(e^{i\th})=U(\cos\th)$, $z_j=e^{i\th_j}$, $\lb_j=\cos\th_j$, $j=0,\dots,r+1$,
and the functions $b_\pm(z)$ are defined in (\ref{WienH}).
\end{theorem}

\begin{remark}
$D_n(1)$ is an explicitly computable determinant related to the Legendre polynomials
(it can also be written as a Selberg integral), c.f. \ci{Warc},
\be
D_n(1)=2^{n^2}\prod_{k=0}^{n-1}
\frac{k!^3}{(n+k)!}
=\frac{\pi^{n+1/2}G(1/2)^2}{2^{n(n-1)}n^{1/4}}\left(1+o(1)\right).
\ee
\end{remark}

\begin{remark}
Since $\bt_j$ enter the symbol only via $e^{\pm i\pi\bt_j}$, the theorem describes
the general case with the exception of the situation when some 
$\Re\bt_j=1/2\mod 1$.
\end{remark}

To prove Theorem \ref{asHankel} we use the fact that
$w(x)$ can be generated by a particular class of functions 
$f(z)$ given by (\ref{fFH}). Namely, we can find an {\it even} function $f$ of $\th$ 
($f(e^{i\th})=f(e^{-i\th})$, $\th\in [0,2\pi)$) 
such that
\be\la{wf}
w(x)={f(e^{i\th})\over|\sin\th|},\qquad x=\cos\th,\quad x\in[-1,1].
\ee
It turns out that we must have $m=2r+1$, $\th_0=0$, 
$\th_{r+1}=\pi$,
$\th_{m+1-j}=2\pi-\th_j$, $j=1,\dots r$. If we denote the beta-parameters of
$f(z)$ by $\wt\bt_j$, we obtain $\wt\bt_0=\wt\bt_{r+1}=0$,
$\wt\bt_j=-\wt\bt_{m+1-j}=-\bt_j$, $j=1,\dots,r$.  
In particular, $\sum_{j=0}^m\wt\bt_j=0$ as remarked above.

We obtain  Theorem \ref{asHankel} from Theorem \ref{asTop} and 
asymptotics for the orthogonal polynomials on the unit circle with weight $f(z)$
using the following connection we establish between Hankel and Toeplitz determinants:
\be\la{131}
D_n(w(x))^2={\pi^{2n}\over 4^{(n-1)^2}}
{(\chi_{2n}+\phi_{2n}(0))^2\over\phi_{2n}(1)\phi_{2n}(-1)}D_{2n}(f(z)),
\qquad n=1,2,\dots,
\ee
where $w(x)$ and $f(z)$ are related by (\ref{wf}).

\begin{remark}
Asymptotics for a subset of symbols (\ref{wFH}) which satisfy a symmetry condition
and have a certain behaviour at the end-points $\pm 1$ were found 
by Basor and Ehrhardt in \ci{BE1}. They use relations between 
Hankel and Toeplitz determinants
which are less general than (\ref{131}) but do not involve polynomials.
For some other related results, see \ci{Geronimo,KVA}.
\end{remark}

\begin{remark}\la{b12}
Asymptotics of a Hankel determinant when some (or all) of $\bt_j$ have the real part 
$1/2$ can be easily obtained. For the corresponding $f(z)$ this implies that
certain $\Re\wt\bt_j=-1/2$ and $\Re\wt\bt_{m+1-j}=1/2$ and the rest 
$\Re\wt\bt_k\in(-1/2,1/2)$. Thus, Theorem \ref{BT} can be used to estimate 
$D_{2n}(f(z))$. For the asymptotics of $\phi_{2n}(z)$ in this case we need
an additional ``correction'' term which is now 
$O(n^{-2\wt\bt_j-1})=O(1)$.   
\end{remark}

\begin{remark}
One can obtain the asymptotics of the polynomials orthogonal 
on the interval $[-1,1]$ with weight (\ref{wFH}) by using our
results for the polynomials $\phi_k(z)$ orthogonal with the
corresponding even weight on the unit circle
and a Szeg\H o relation which maps 
the latter polynomials to the former ones.
\end{remark}

Our final task is to present asymptotics for the so-called Toeplitz+Hankel 
determinants.
We consider the four most important ones appearing in the theory of classical groups and its applications to random matrices and statistical mechanics
(see, e.g., \ci{Baik,FF,KM}) defined in terms of
the Fourier coefficients of an even $f$ (evenness implies
the matrices are symmetric) as follows:
\be\la{intT+H}
\det (f_{j-k}+f_{j+k})_{j,k=0}^{n-1},\quad
\det (f_{j-k}-f_{j+k+2})_{j,k=0}^{n-1},\quad
\det (f_{j-k}\pm f_{j+k+1})_{j,k=0}^{n-1}.
\ee
There are simple relations \ci{Weyl,J2,Baik} between the determinants (\ref{intT+H}) 
and Hankel determinants on $[-1,1]$ with added singularities at the end-points,
namely,
\begin{eqnarray}
\det (f_{j-k}+f_{j+k})_{j,k=0}^{n-1}={2^{n^2-2n+2}\over\pi^n} 
D_n(f(e^{i\theta(x)})/\sqrt{1-x^2}),\\
\det (f_{j-k}-f_{j+k+2})_{j,k=0}^{n-1}={2^{n^2}\over\pi^n} 
D_n(f(e^{i\theta(x)})\sqrt{1-x^2}),\la{TH-}\\
\det (f_{j-k}+f_{j+k+1})_{j,k=0}^{n-1}={2^{n^2-n}\over\pi^n} 
D_n\left(f(e^{i\theta(x)})\sqrt{1+x\over 1-x}\right),\la{TH+1}\\
\det (f_{j-k}-f_{j+k+1})_{j,k=0}^{n-1}={2^{n^2-n}\over\pi^n} 
D_n\left(f(e^{i\theta(x)})\sqrt{1-x\over 1+x}\right),\la{TH-1}
\end{eqnarray}
where $f$ is even, and $x=\cos\theta$.
It is easily seen that if $f(z)$ is an (even) function of type (\ref{fFH})
then the corresponding symbols of Hankel determinants belong to 
the class (\ref{wFH}). Thus a combination of the above formulas and
Theorem \ref{asHankel} gives the following

\begin{theorem}\la{T+H} 
Let $f(z)$ be defined in (\ref{fFH}) with the condition
$f(e^{i\th})=f(e^{-i\th})$. Let $\th_{r+1}=\pi$. Then as $n\to\infty$,
\begin{multline}\la{TH}
D_n^{\mathrm{T+H}}=
e^{nV_0+{1\over 2}\left[(\al_0+\al_{r+1}+s+t)V_0-(\al_0+s)V(1)-(\al_{r+1}+t)V(-1)+
\sum_{k=1}^\infty k V_k^2\right]}\\
\times
\prod_{j=1}^r b_+(z_j)^{-\al_j+\bt_j}b_-(z_j)^{-\al_j-\bt_j}\times
e^{-i\pi\left[\left\{\al_0+s+\sum_{j=1}^r\al_j\right\}\sum_{j=1}^r\bt_j+
\sum_{1\le j<k\le r}(\al_j\bt_k-\al_k\bt_j)
\right]}\\
\times
2^{(1-s-t)n+p+\sum_{j=1}^r(\al_j^2-\bt_j^2)-{1\over2}(\al_0+\al_{r+1}+s+t)^2
+{1\over2}(\al_0+\al_{r+1}+s+t)}
n^{{1\over2}(\al_0^2+\al_{r+1}^2)+\al_0 s+\al_{r+1}t+
\sum_{j=1}^r(\al_j^2-\bt_j^2)}\\
\times
\prod_{1\le j<k\le r}|z_j-z_k|^{-2(\al_j\al_k-\bt_j\bt_k)}
|z_j-z_k^{-1}|^{-2(\al_j\al_k+\bt_j\bt_k)}\\
\times
\prod_{j=1}^r z_j^{2\wt A\bt_j}|1-z_j^2|^{-(\al_j^2+\bt_j^2)}
|1-z_j|^{-2\al_j(\al_0+s)}|1+z_j|^{-2\al_j(\al_{r+1}+t)}\\
\times
\frac{\pi^{{1\over2}(\al_0+\al_{r+1}+s+t+1)}G(1/2)^2}
{G(1+\al_0+s)G(1+\al_{r+1}+t)}
\prod_{j=1}^r\frac{G(1+\al_j+\bt_j) G(1+\al_j-\bt_j)}{G(1+2\al_j)}
\left(1+o(1)\right),\\
\wt A={1\over2}(\al_0+\al_{r+1}+s+t)+\sum_{j=1}^r\al_j,\\
\Re\al_j>-{1\over 2},\qquad \Re\bt_j\in\left(-{1\over 2},{1\over2}\right),
\qquad j=0,1,\dots,r+1,\qquad \bt_0=\bt_{r+1}=0.
\end{multline}
Here
\begin{align}
&D_n^{\mathrm{T+H}}=\det (f_{j-k}+f_{j+k})_{j,k=0}^{n-1},\qquad
\mbox{with}\quad p=-2n+2,\quad s=t=-{1\over 2}\\
&D_n^{\mathrm{T+H}}=\det (f_{j-k}-f_{j+k+2})_{j,k=0}^{n-1},\qquad
\mbox{with}\quad p=0,\quad s=t={1\over 2}\\
&D_n^{\mathrm{T+H}}=\det (f_{j-k}\pm f_{j+k+1})_{j,k=0}^{n-1},\qquad
\mbox{with}\quad p=-n,\quad s=\mp{1\over 2},\quad t=\pm{1\over 2}.
\end{align}
\end{theorem}

\begin{remark}
For the case $\Re\bt_j=1/2$ see Remark \ref{b12} above.
\end{remark}

\begin{remark}
For the determinant $\det (f_{j-k}+f_{j+k+1})_{j,k=0}^{n-1}$ in the case
when the symbol has no $\al$ singularities at $z=\pm 1$ and $|\Re\bt_j|<1/2$, 
the asymptotics 
were obtained in \ci{BEsym} (see also \ci{BEnonsym} if $f$ is non-even, $\al_j=0$).
Note that for symbols without singularities, i.e. for $f(z)=e^{V(z)}$,
the asymptotics of all the above Toeplitz+Hankel determinants 
(and related more general ones) were found recently in \ci{BE4}.
\end{remark}


\section*{Acknowledgements}
Percy Deift was supported
in part by NSF grant \# DMS 0500923. 
Alexander Its was supported
in part by NSF grant \# DMS-0701768 and
EPSRC grant \# EP/F014198/1. Igor Krasovsky was 
supported in part by EPSRC grants \# EP/E022928/1
and \# EP/F014198/1.

\end{document}